\documentclass[12pt,a4paper]{amsart}
\usepackage{latexsym,amsmath,amsfonts,amscd,amssymb}
\theoremstyle{plain}  
\newtheorem{theorem}{Theorem}[section]

\newtheorem*{theorem*}{Theorem}

\newtheorem{lemma}[theorem]{Lemma}
\newtheorem{proposition}[theorem]{Proposition}

\theoremstyle{remark}

\newtheorem{remark}[theorem]{Remark}

\newtheorem*{claim*}{Claim}
\numberwithin{equation}{section}
\newcommand{\suchthat}{\;\;|\;\;}

\renewcommand{\leq}{\leqslant}
\renewcommand{\le}{\leqslant}
\renewcommand{\geq}{\geqslant}

\renewcommand{\setminus}{\smallsetminus}

\newcommand{\lto}{\longrightarrow}
\newcommand{\R}{\mathbb{R}}

\newcommand{\Z}{\mathbb{Z}}
\newcommand{\C}{\mathbb{C}}
\newcommand{\Q}{\mathbb{Q}}
\newcommand{\HH}{\mathbb{H}}
\newcommand{\PP}{\mathbb{P}}

\newcommand{\dbar}{\bar{\partial}}

\newcommand{\lra}{\longrightarrow}
\newcommand{\PU}{\mathrm{PU}}
\newcommand{\PGL}{\mathrm{PGL}}

\newcommand{\SU}{\mathrm{SU}}
\newcommand{\U}{\mathrm{U}}
\newcommand{\GL}{\mathrm{GL}}
\newcommand{\SL}{\mathrm{SL}}

\newcommand{\Sp}{\mathrm{Sp}}
\newcommand{\GCD}{\mathrm{GCD}}

\DeclareMathOperator{\ad}{ad}

\DeclareMathOperator{\tr}{tr}
\DeclareMathOperator{\rk}{rk}
\DeclareMathOperator{\im}{im}
\DeclareMathOperator{\coker}{coker}
\DeclareMathOperator{\Hom}{Hom}
\DeclareMathOperator{\End}{End}

\DeclareMathOperator{\codim}{codim}
\DeclareMathOperator{\Index}{index}

\newcommand{\noi}{\noindent}

\newcommand{\norm}[1]{\lVert#1\rVert}

\newcommand{\Gr}{\operatorname{Gr}}

\begin{document}
\title{Homotopy groups of moduli spaces of representations}

\date{10/07/2006}

\author[S. B. Bradlow]{Steven B. Bradlow}
\address{Department of Mathematics \\
University of Illinois \\
Urbana \\
IL 61801 \\
USA }
\email{bradlow@math.uiuc.edu}

\author[O. Garc{\'\i}a-Prada]{Oscar Garc{\'\i}a-Prada}
\address{Instituto de Matem\'aticas \\
  CSIC \\ Serrano 121 \\ 28006 Madrid \\ Spain}
\email{oscar.garcia-prada@uam.es}

\author[P. B. Gothen]{Peter B. Gothen}
\address{Departamento de Matem\'atica Pura \\
  Faculdade de Ci\^encias da Universidade do Porto \\
  Rua do Campo Alegre 687 \\ 4169-007 Porto \\ Portugal }
\email{pbgothen@fc.up.pt}

\thanks{
  Members of VBAC (Vector Bundles on Algebraic Curves).
  Second and Third authors partially supported by Ministerio de
  Educaci\'{o}n y Ciencia and Conselho de Reitores das
  Universidades Portuguesas through Acci\'{o}n Integrada Hispano-Lusa
  HP-2000-0015 (Spain) / E--30/03 (Portugal).
  Second author partially supported by Ministerio de Educaci{\'o}n y
  Ciencia (Spain) through Project BFM2000-0024.
  Third author partially supported by the Centro de Matem\'atica da
  Universidade do Porto and the project POCTI/MAT/58549/2004, financed
  by FCT (Portugal) through the programmes POCTI and POSI of the QCA
  III (2000--2006) with European Community (FEDER) and national funds.
}



\begin{abstract}
  We calculate certain homotopy groups of the moduli spaces for
  representations of a compact oriented surface in the Lie groups
  $\GL(n,\C)$ and $\U(p,q)$.  Our approach relies on the
  interpretation of these representations in terms of Higgs bundles
  and uses Bott--Morse theory on the corresponding moduli spaces.
\end{abstract}

\maketitle

\section{Introduction}

Given a closed oriented surface, $X$, and a Lie group $G$, moduli spaces of
surface group representations in $G$ have rich geometric and topological
structure which reflects properties of both $X$ and $G$. In this paper we
consider the cases where $G$ is $\GL(n,\C)$ or $\U(p,q)$.

Our main tools rely on an interpretation of the moduli spaces in terms of
holomorphic bundles. Such an interpretation starts from the basic
correspondence between representations of the fundamental group and flat
principal bundles. Holomorphic bundles enter the picture if we fix a
complex structure on the surface $X$ --- thereby turning it into a Riemann
surface. By results of Hitchin \cite{H}, Donaldson \cite{D}, Simpson
\cite{S1} and Corlette \cite{C}, if $G$ is complex semisimple then the
flat principal $G$-bundles corresponding to semisimple representations
of $\pi_1 X$
in $G$ are equivalent to polystable $G$-Higgs bundles over the Riemann
surface. More generally, such Higgs bundles exist if $G$ is complex
reductive, in which case the polystable $G$-Higgs bundles correspond to
semisimple representations not of $\pi_1 X$ but of a central extension of
the fundamental group.

Referring to $\pi_1 X$ and its central extensions as surface groups, we can
thus identify the moduli spaces of surface group representations with
moduli spaces of polystable Higgs bundles. This identification puts a
natural K\"ahler structure on the moduli spaces and also reveals a
compatible $\C^*$-action. The restriction of this action to $S^1$ leads to
a symplectic moment map whose squared norm serves as a proper Morse
function. In a striking example of the interplay between geometry and
topology, these geometric features on the moduli space of Higgs bundles
provide powerful tools for studying the topology of the underlying moduli
spaces of surface group representations.

Holomorphic bundle techniques can also be adapted to the case in which
$G$ is a real reductive Lie group, in particular when $G$ is a real
form of a complex reductive group. If $G$ is the compact real form of
a complex reductive group $G_{\C}$, then the theorem of Narasimhan and
Seshadri \cite{NS} and its generalization by Ramanathan \cite{R}
identify representations into $G$ with polystable principal
$G_{\C}$-bundles.  For non-compact real forms the basic ideas were
first introduced by Nigel Hitchin. In \cite{H2} he outlined how to
define the appropriate Higgs bundles and applied his methods to the
case $G=\SL(n,\R)$ and also to other split real forms. Other special
cases have been considered in a similar
way\footnotemark\footnotetext{Notably $\Sp(4,\R)$ and $\SU(2,2)$
  \cite{gothen:2001, garcia-prada-mundet:2004}, $\U(2,1)$
  \cite{gothen:2002}, $\U(p,q)$ and $\PU(p,q)$
  \cite{BGG1},
  $\GL(n,\R)$ \cite{bradlow-garcia-prada-gothen:2004b}. Higgs bundle
  methods have also been applied, albeit in a more algebraic way in
  the cases $\U(p,1)$ \cite{xia:1999}, $\PU(2,1)$ \cite{xia:2000}, and
  $\PU(p,p)$ \cite{markman-xia:2002}.}. In \cite{BGG1} we began an in-depth
study of the groups $\U(p,q)$ and their adjoint forms $\PU(p,q)$ (for
any $p$ and $q$) from this point of view. This paper is a continuation
of that work.

The most primitive topological feature of the moduli spaces is their number
of connected components, i.e.\ $\pi_0$. The above methods have been
effective in addressing this question, mainly by exploiting the properness
of the above mentioned Morse function. This transfers questions about
$\pi_0$ for the moduli spaces into questions about the connected components
of the minimal submanifolds for the Morse function.

In good cases, there is additional useful Morse theoretic information which
has thus far gone unexploited. Our goal is to correct this oversight. In
particular, using information about the Morse indices of non-minimal
critical points,  we can relate higher homotopy groups for the full moduli
spaces to those of their minimal submanifolds. For the latter we rely on
the calculations by Daskalopoulos and Uhlenbeck \cite{DU} for higher
homotopy groups of the moduli space of stable vector bundles. Our main
results for the moduli spaces of $\GL(n,\C)$ and $\U(p,q)$ Higgs bundles,
and hence for the corresponding moduli spaces of representations, are given
in Theorems \ref{homotopy-higgs-n-d} and \ref{homotopy-higgs-p-q-a-b}
respectively.

\textbf{Acknowledgments.} We thank Tam\'as Hausel for enlightening
discussions.  We are greatly indebted to an anonymous referee for
pointing out that our initial estimate in (\ref{item:index-gln}) of
Proposition~\ref{index-bound} could be improved and for providing the
argument outlined in Remark~\ref{rem:referee}.

\section{Surface group representations and Higgs bundles}\label{sec:2}

For a more thorough account of the material in this section see
\cite{BGG1}.

\subsection{Surface group representations}

Let $X$ be a smooth closed oriented surface of genus $g\geq 2$. The
fundamental group, $\pi_1 X$, of $X$ is a finitely generated group
generated by $2g$ generators, say $A_{1},B_{1}, \ldots, A_{g},B_{g}$,
subject to the single relation $\prod_{i=1}^{g}[A_{i},B_{i}] = 1$. It has a
universal central extension
\begin{equation}\label{eq:gamma}
0\longrightarrow\mathbb{Z}\longrightarrow\Gamma\longrightarrow\pi_1
X\longrightarrow 1 \
\end{equation}

\noindent generated by the same generators as $\pi_1 X$, together with a
central element $J$ subject to the relation $\prod_{i=1}^{g}[A_{i},B_{i}] =
J$.

By a representation of $\Gamma$ in $\GL(n,\C)$ we mean a homomorphism $\rho
\colon \Gamma\to \GL(n,\C)$.  We say that a representation of $\Gamma$ in
$\GL(n,\C)$ is \emph{semisimple} if the $\C^n$-representation of
$\Gamma$ induced by the fundamental representation of
$\GL(n,\C)$ is semisimple\footnote{In general a representation of
  $\Gamma$ in a reductive Lie group $G$ is said to be semisimple if the induced
  (adjoint) representation on the Lie algebra of $G$ is
  semisimple. For $G=\GL(n,\C)$ this is equivalent to the
  definition given here.}.  The group $\GL(n,\C)$
acts on the set of representations via conjugation. Restricting to the
semisimple representations, denoted by $\Hom^{+}(\Gamma, \GL(n,\C))$, we
get the \emph{moduli space} of representations of $\Gamma$ in $\GL(n,\C)$,

\begin{equation}\label{eqn:RGdef}
  \mathcal{R}(\Gamma,\GL(n,\C)) = \Hom^+(\Gamma, \GL(n,\C)) / \GL(n,\C).
\end{equation}

\noi The set $\Hom^{+}(\Gamma, \GL(n,\C))$ can be embedded in
$\GL(n,\C)^{2g+1}$ via the map
\begin{align*}
  \Hom^{+}(\Gamma, \GL(n,\C)) &\to \GL(n,\C)^{2g+1} \\
  \rho &\mapsto (\rho(A_1), \ldots \rho(B_g), \rho(J)).
\end{align*}

\noindent We  can then give $\Hom^{+}(\Gamma, \GL(n,\C))$ the subspace
topology and $\mathcal{R}(\Gamma,\GL(n,\C))$ the quotient topology. This
topology is Hausdorff because we have restricted attention to semisimple
representations.

There is a topological invariant of a representation $\rho\in
\mathcal{R}(\Gamma,\GL(n,\C))$ given  by $\rho(J)$,  which coincides with
the first Chern class of the vector bundle with central curvature
associated to $\rho$. Fixing this invariant, we define
$$
\mathcal{R}(n,d):=\{ \rho\in \mathcal{R}(\Gamma, \GL(n,\C))\;\;|\;\;
\rho(J) = [d] \in
\Z_n\subset  Z(\GL(n,\C)\}.
$$

\noi In particular the representations with vanishing degree correspond to
representations of the fundamental group of $X$, that is,
\begin{equation}\label{eqn:RGdef2}
\mathcal{R}(n,0)= \mathcal{R}(\pi_1 X, \GL(n,\C)) := \Hom^+(\pi_1 X,
\GL(n,\C)) / \GL(n,\C)\ .
\end{equation}

\bigskip

Similarly to the case of $\GL(n,\C)$ we consider the moduli space
\begin{equation}\label{eqn:RGdef3}
  \mathcal{R}(\Gamma, \U(p,q)) = \Hom^+(\Gamma, \U(p,q)) / \U(p,q)\ .
\end{equation}

\noi The moduli space $\mathcal{R}(\Gamma,\U(p,q))$ can be identified with
the moduli space of $\U(p,q)$-bundles on $X$ with projectively flat
connections.  Taking a reduction to the maximal compact $\U(p)\times\U(q)$,
we thus associate to each class $\rho\in \mathcal{R}(\Gamma, \U(p,q))$ a
vector bundle of the form $V\oplus W$, where $V$ and $W$ are rank $p$ and
$q$ respectively, and thus a pair of integers $(a,b) =(\deg(V),\deg(W))$.
There is thus a map
\begin{displaymath}
  {c} \colon  \mathcal{R}(\Gamma,\U(p,q)) \to \Z\oplus\Z
\end{displaymath}
given by ${c}(\rho) = (a,b)$. The corresponding map on $\Hom^+(\Gamma,
\U(p,q))$ is clearly continuous and thus locally constant. Since $\U(p,q)$
is connected, the map ${c}$ is likewise continuous and thus constant on
connected components. The subspace of $\mathcal{R}(\Gamma,(\U(p,q))$
corresponding to representations with invariants $(a,b)$ is denoted by
\begin{equation}\label{R(p,q,a,b)}
  \mathcal{R}(p,q, a,b) = {c}^{-1}(a,b)=
  \{ {\rho} \in \mathcal{R}(\Gamma,\U(p,q)) \suchthat
  {c}({\rho})=(a,b) \in \Z \oplus \Z\}\ .
\end{equation}

\noindent The representations for which $a+b=0$  correspond to
representations of the fundamental group of $X$, that is,
\begin{equation}\label{eqn:RGdef4}
\ \mathcal{R}(p,q, a,-a) = {c}^{-1}(a,-a)=
  \{ {\rho} \in \mathcal{R}(\pi_1 X ,\U(p,q)) \suchthat
  {c}({\rho})=(a,-a) \in \Z \oplus \Z\}\ .
\end{equation}

\subsection{$\GL(n,\C)$-Higgs bundles}

A $\GL(n,\C)$-Higgs bundle on a compact Riemann surface $X$ is a pair
$(E,\Phi)$, where $E$ is a rank $n$ holomorphic vector bundle over $X$ and
$\Phi \in H^0(\End(E) \otimes K)$ is a holomorphic endomorphism of $E$
twisted by the canonical bundle $K$ of $X$. The $\GL(n,\C)$-Higgs bundle
$(E,\Phi)$ is \emph{stable} if
  the slope stability condition
\begin{equation}\label{eq:stability}
\mu(E') < \mu(E)
\end{equation}
holds for all proper $\Phi$-invariant subbundles $E'$ of $E$.  Here
the \emph{slope} is defined by $\mu(E)=\deg(E)/\rk(E)$ and
\emph{$\Phi$-invariance} means that $\Phi(E')\subset E'\otimes K$.
\emph{Semistability} is defined by replacing the above strict inequality
with a weak inequality. A Higgs bundle is called \emph{polystable} if it is
the direct sum of stable Higgs bundles with the same slope.

Given a hermitian metric on $E$, let $A$ denote the unique unitary
connection compatible with the holomorphic structure, and let
  $F_A$ be its curvature.
  \emph{Hitchin's equations} on $(E,\Phi)$ are
\begin{equation}
  \label{eq:hitchin1}
  \begin{aligned}
    F_A + [\Phi,\Phi^*] &= -\sqrt{-1}\mu \text{Id}_E \omega, \\
    \dbar_{A} \Phi &=0,
  \end{aligned}
\end{equation}
where $\omega$ is the K\"ahler form on $X$, $\text{Id}_E$ is the identity
on $E$, $\mu = \mu(E)$ and $\dbar_A$ is the anti-holomorphic part of the
covariant derivative $d_A$.  A solution to Hitchin's equations is
\emph{irreducible} if there is no proper subbundle of $E$ preserved by $A$
and $\Phi$.

If we define a Higgs connection (as in
  \cite{S2}) by
\begin{equation}\label{eqn:Dhiggs}
D=d_A+\theta
\end{equation}
\noindent where $\theta=\Phi+\Phi^* $, then Hitchin's equations are
equivalent to the conditions
\begin{equation}
  \label{eq:harmonic1}
  \begin{aligned}
  F_D &= -\sqrt{-1}\mu \text{Id}_E \omega,\\
  d_A^* \theta &=0.
  \end{aligned}
\end{equation}
In particular, the first equation says that $D$ is a projectively flat
connection\footnotemark\footnotetext{ The other equation is an
harmonicity constraint.}. If $\deg(E)=0$ then $D$ is actually flat. It
follows that in this case the pair $(E,D)$ defines a representation of
$\pi_1 X$ in $\GL(n,\C)$. If $\deg(E)\ne 0$, then the pair $(E,D)$ defines
a representation of $\pi_1 X$ in $\PGL(n,\C)$, or equivalently, a
representation of $\Gamma$ in $\GL(n,\C)$. By the theorem of Corlette
(\cite{C}), every semisimple representation of $\Gamma$ (and therefore
every semisimple representation of $\pi_1 X$) arises in this way.

If we fix the rank and degree (say $n$ and $d$ respectively) of the bundle
$E$, i.e.\ on bundles of fixed topological type, the isomorphism classes of
polystable Higgs bundles are parameterized by a quasi-projective variety of
dimension $2 +2 n^2(g-1)$. We denote this moduli space of rank $n$ degree
$d$ polystable Higgs bundles by $\mathcal{M}(n,d)$.

If we fix a hermitian metric on a smooth rank $n$ degree $d$ complex vector
bundle on $X$, then there is a gauge theoretic moduli space of pairs
$(A,\Phi)$, consisting of a unitary connection $A$ and an endomorphism
valued $(1,0)$-form $\Phi$, which are solutions to Hitchin's equations
\eqref{eq:hitchin1}, modulo $\U(n)$-gauge equivalence.

The gauge theory moduli space and $\mathcal{M}(n,d)$ are related by virtue
of the Hitchin-Kobayashi correspondence: a $\GL(n,\C)$-Higgs bundle
$(E,\Phi)$ is polystable if and only if it admits a hermitian metric such
that Hitchin's equations \eqref{eq:hitchin1} are satisfied, and $(E,\Phi)$
is stable if and only if the corresponding solution is irreducible. There
is, moreover, a map from the gauge theoretic moduli space to this moduli
 space given by taking a solution $(A,\Phi)$ to Hitchin's equations to the
 Higgs bundle $(E,\Phi)$, where the holomorphic structure on $E$ is
 given by $\dbar_{A}$.  This map is a homeomorphism, and a diffeomorphism
 on the smooth locus.

 In view of the relation between Hitchin's equations and projectively
 flat connections, this correspondence gives rise to a homeomorphism
 between $\mathcal{M}(n,d)$ and the component $\mathcal{R}(n,d)$ of
 the moduli space of semisimple representations of $\Gamma$ in
 $\GL(n,\C)$.  If the degree of the Higgs bundle is zero, then the
 moduli space $\mathcal{M}(n,0)$ is homeomorphic to the moduli space
 of representations of $\pi_1 X$ in $\GL(n,\C)$.

 \begin{theorem} If $(n,d)$ is such that $\GCD(n,d)=1$ then the moduli
   space $\mathcal{M}(n,d)$ is a non-empty connected smooth
   hyperk\"ahler manifold.
\end{theorem}

\subsection{$\U(p,q)$-Higgs bundles}\label{subs:UpqHiggs}

There is a special class   of $\GL(n,\C)$-Higgs bundles, related to
representations in $\U(p,q)$ given  by the requirements that
\begin{equation} \label{upq-higgs-bundle}
  \begin{aligned}
  E &= V \oplus W \\
  \Phi &=
  \left(
  \begin{matrix}
    0 & \beta \\
    \gamma & 0
  \end{matrix}
  \right)
  \end{aligned}
\end{equation}
where $V$ and $W$ are holomorphic vector bundles of rank $p$ and $q$
respectively and the non-zero components in the Higgs field are $\beta \in
H^0(\Hom(W,V)\otimes K)$, and $\gamma \in H^0(\Hom(V,W) \otimes K)$. We say
$(E,\Phi)$ is a \emph{stable} $\U(p,q)$-Higgs bundle if the slope stability
condition $\mu(E') < \mu(E)$, is satisfied for all $\Phi$-invariant
subbundles $E'=V'\oplus W'$, i.e.\ for all subbundles $V'\subset V$ and
$W'\subset W$ such that
\begin{align}
\beta &:W'\longrightarrow V'\otimes K\\
\gamma &:V'\longrightarrow W'\otimes K\ .
\end{align}

\noi Semistability and polystability are defined analogously to the way
they are defined for $\GL(n,\C)$-Higgs bundles.

Hitchin's equations make sense for $\U(p,q)$-Higgs bundles, with a $\U(p,q)$
solution being a metric with respect to which $E= V \oplus W$ is an
orthogonal decomposition. With $\Phi$ as in (\ref{upq-higgs-bundle}) and
$\theta=\Phi+\Phi^* $, the corresponding $\U(p,q)$-Higgs connection
$D=d_A+\theta$ is not only projectively flat but has $\U(p,q)$ holonomy.
This provides the link between $\U(p,q)$-Higgs bundles and surface group
representations in $\U(p,q)$, leading to:

\begin{theorem}\label{prop:R=M} Let $\mathcal{M}(p,q,a,b)$ be the moduli
space of polystable $\U(p,q)$-Higgs bundles with $\deg(V)=a$ and $\deg
W=b$. Then with $\mathcal{R}(p,q,a,b)$ as in (\ref{R(p,q,a,b)}) there is a
homeomorphism $\mathcal{M}(p,q,a,b)\cong \mathcal{R}(p,q,a,b)$.
 \end{theorem}

The \emph{Toledo invariant} of the representation corresponding
to $(E=V\oplus W,\Phi)$ is defined by
  \begin{equation}
\tau=\tau(p,q,a,b)    = 2\frac{qa - pb}{p + q}
\label{toledo-invariant}
  \end{equation}
\noindent where $a=\deg(V)$ and $b=\deg(W)$. This invariant satisfies the
following  Milnor-Wood-type inequality (proved by Domic and Toledo
\cite{DT})

\begin{equation}
|\tau(p,q,a,b)|\leq \min\{p,q\}(2g-2)\ .
\end{equation}

\begin{theorem}\cite{BGG1}
Let $(p,q,a,b)$ such that $\GCD(p+q,a+b)=1$. Then $\mathcal{M}(p,q,a,b)$
(and hence $\mathcal{R}(p,q,a,b)$) is a connected smooth K\" ahler manifold
which is non-empty if and only if $|\tau(p,q,a,b)|\leq \min\{p,q\}(2g-2)$.
\end{theorem}

\section{Morse theory on the moduli space}

\subsection{The Morse function}
\label{sec:morse-function}

Let ${\mathcal M}$ be either $\mathcal{M}(n,d)$ or ${\mathcal M}(p,q,a,b)$.
We will assume that $\GCD(n,d) = 1$ and $\GCD(p+q,a+b)=1$. Under this
coprimality condition, there are no strictly semistable Higgs bundles and
the moduli space ${\mathcal M}$ is smooth. The non-zero complex numbers
$\C^*$ act on $\mathcal{M}$ via the map $\lambda\cdot(E,\Phi) = (E, \lambda
\Phi)$.  However, to have an action on the gauge theory moduli space (i.e.\
on the set of solutions to Hitchin's equations \eqref{eq:hitchin1}, cf.\
Section~\ref{sec:2}), one must restrict to the action of $S^1 \subset
\C^*$.  This is a Hamiltonian action and the associated moment map is
\begin{displaymath}
  [(A,\Phi)] \mapsto -\tfrac{1}{2}\norm{\Phi}^2
  = -i \int_X \tr (\Phi\Phi^*)
\end{displaymath}
where the adjoint $\Phi^*$ is taken with respect to the hermitian
metric on $E$.
We shall, however, prefer to consider the positive function
\begin{equation}\label{eq:aaa4}
  f([A,\Phi]) = \tfrac{1}{2}\norm{\Phi}^2.
\end{equation}

Next we recall a general result of Frankel \cite{F}, which was first
used in the context of moduli spaces of Higgs bundles by Hitchin
\cite{H}.

\begin{theorem}
  Let $\tilde{f}\colon M \to \R$ be a proper moment map for a
  Hamiltonian circle action on a K\"ahler manifold $M$.  Then
  $\tilde{f}$ is a perfect Bott--Morse function.  
\end{theorem}

\subsection{Morse theory and homotopy groups}
\label{sec:morse-homotopy}

In this Section we recall some basic facts of Bott--Morse theory.
Let ${\mathcal M}_l\subset \mathcal{M}$ be the critical submanifolds
of $f$ and $\nu(\mathcal{M}_l)$ be the normal bundle of ${\mathcal
  M}_l$ in $\mathcal{M}$.  The Hessian of $f$ is non-degenerate on
$\nu(\mathcal{M}_l)$ and we have the decomposition in positive and
negative eigenspace bundles
$$
\nu(\mathcal{M}_l)= \nu^+(\mathcal{M}_l)\oplus  \nu^-(\mathcal{M}_l).
$$
The index of $\mathcal{M}_l$ is defined as
$$
\Index(\mathcal{M}_l):=\rk \nu^-(\mathcal{M}_l).
$$
Let ${\mathcal M}_l^+$ be the stable set of $\mathcal{M}_l$, i.e.,
the subset of $\mathcal{M}$ defined by the points of $\mathcal M$
which flow to $\mathcal M_l$.  It follows from Bott--Morse theory that
${\mathcal M}_l^+$ is a submanifold of $\mathcal{M}$ of codimension
\begin{equation}
  \label{eq:codimension-index}
  \codim_\R(\mathcal{M}_l^+)= \Index(\mathcal{M}_l),
\end{equation}
and that there is a stratification
\begin{equation}\label{eq:morse-stratification}
  \mathcal{M} = \bigcup_l {\mathcal M_l^+}.
\end{equation}

\begin{proposition}\label{homotopy-minima}
  Let $\mathcal N = \mathcal{M}_0 \subset \mathcal{M}$ be the
  submanifold of local minima of $f$. If $\Index(\mathcal{M}_l)\geq
  m\geq 2$ for every $l \neq 0$ then
   $$
     \pi_i(\mathcal{M})\cong \pi_i(\mathcal{N})
     \;\;\;\text{for}\;\;\; i\leq m-2.
   $$
\end{proposition}
\begin{proof}
  The stratification \eqref{eq:morse-stratification} shows that
  $$
    \mathcal{M}_0^+={\mathcal M} \setminus \bigcup_{l\neq 0} {\mathcal M_l^+}
  $$
  and the Morse flow defines a retraction from $\mathcal M_0^+$ to
  $\mathcal{N} = \mathcal{M}_0$.  Thus the result is an immediate
  consequence of standard homotopy theory, using
  \eqref{eq:codimension-index}.
\end{proof}

\subsection{Deformation theory of Higgs bundles}
\label{sec:deformation}

In the following we recall some standard facts about the deformation
theory of Higgs bundles (this has been treated in many places, a
convenient reference is Biswas--Ramanan \cite{biswas-ramanan:1994}).
In order to describe the results in a uniform way for a $G$-Higgs bundle
$(E,\Phi)$ when $G=\GL(n,\C)$ or $\U(p,q)$, we introduce bundles
$U_G^+$, $U_G^-$ and $U_G$ defined by
\begin{align*}
  U_{\GL(n,\C)}^+ &= U_{\GL(n,\C)}^- = U_{\GL(n,\C)}  = \End(E), \\
  U_{\U(p,q)}^+ &= \End(V) \oplus \End(W), \\
  U_{\U(p,q)}^- &= \Hom(W,V) \oplus \Hom(V,W), \\
  U_{\U(p,q)} &= U_{\U(p,q)}^+ \oplus U_{\U(p,q)}^- =\End(V \oplus W),
\end{align*}

\noi where the bundles $V$ and $W$ are as in Section \ref{subs:UpqHiggs}.
 Note that, with this notation, $\Phi \in H^0(U^-_G\otimes K)$.

\begin{remark}\label{rem:duality}
  Both for $G = \GL(n,\C)$ and for $G = \U(p,q)$, there is an inner
  product on $U_G$ which is invariant under the adjoint action of
  $U_G$, i.e.,
  \begin{equation}\label{eq:inner-invariant}
    \langle \ad(\psi)x,y \rangle + \langle x,\ad(\psi)y \rangle = 0
  \end{equation}
  for local sections $x$, $y$ and $\psi$ of $U_G$.  This inner
  product restricts to an inner product on $U^{-}_G$ and
  $U^{+}_G$, giving rise to an isomorphism
  \begin{equation}\label{eq:duality}
    U^{\pm}_G \xrightarrow{\cong} (U^{\pm}_G)^*.
  \end{equation}
  Note that under this duality
  $$\ad(\Phi)^t = -\ad(\Phi)\otimes 1_{K^{-1}}.$$
\end{remark}

\begin{proposition}\label{prop:deform}
  Let $(E,\Phi)$ be a $G$-Higgs bundle for $G=\GL(n,\C)$ or
  $G=\U(p,q)$ and define the following complex of sheaves
  \begin{displaymath}
    C^{\bullet}_{G}(E,\Phi)\colon
      U_G^+ \xrightarrow{\ad(\Phi)} U_G^- \otimes K.
  \end{displaymath}
  Then the following holds:
  \begin{enumerate}
  \item The space of endomorphisms of $(E,\Phi)$ is naturally
    isomorphic to $\HH^0(C^{\bullet}_{G})$.
  \item The infinitesimal deformation space of $(E,\Phi)$ is naturally
    isomorphic to $\HH^1(C^{\bullet}_{G})$.
  \end{enumerate}
\end{proposition}

The following proposition is simply a statement of the fact that a
stable Higgs bundle is simple.
\begin{proposition}\label{prop:duality}
  Let $(E,\Phi)$ be a stable $G$-Higgs bundle for $G=\GL(n,\C)$ or
  $G=\U(p,q)$.  Then
  \begin{displaymath}
    \HH^0(C^{\bullet}_{G}(E,\Phi)) \cong \C,
  \end{displaymath}
  generated by the scalar multiples of the identity morphism. 
\end{proposition}

\subsection{Critical points and Morse indices}
\label{sec:critical-morse-indices}
In the following $(E,\Phi)$ continues to denote a $G$-Higgs bundle for
$G=\GL(n,\C)$ or $G=\U(p,q)$ and for ease of notation we omit the
subscript $G$ on the bundles $U^\pm_G$ and the complex
$C^{\bullet}_G$.  The critical points of the function $f$ are exactly
the fixed points of the $S^1$-action on $\mathcal{M}$.  This allows
one to describe the corresponding Higgs bundles as ``complex
variations of Hodge structure'', as follows (cf.\ Hitchin \cite{H,H2}
and also Simpson \cite{S2}).

\begin{proposition}
If $(E,\Phi)$ corresponds to a critical
point of $f$, then there is a semisimple element $\psi\in H^0(U^+)$
and a corresponding decomposition in eigenspace bundles
\begin{equation}\label{eq:vhs}
  U^{\pm}_G = \bigoplus_k U^{\pm}_k
\end{equation}
for the adjoint action of $\psi$, such that $\ad(\psi)$ has eigenvalue
$ik$ on $U^\pm_k$.  Furthermore, $[\psi,\Phi]=i\Phi$, i.e.,
\begin{displaymath}
  \Phi \in H^0(U^-_1\otimes K).
\end{displaymath}
\end{proposition}
In particular, this means that the deformation complex of $(E,\Phi)$
decomposes as
\begin{equation}\label{eq:complex-dec}
  C^{\bullet}(E,\Phi) = \bigoplus_k C^{\bullet}_k(E,\Phi),
\end{equation}
where we have defined for each $k$ the complex
\begin{displaymath}
  C^{\bullet}_k(E,\Phi)\colon U_k^+
    \xrightarrow{\ad(\Phi)}U_{k+1}^-\otimes K.
\end{displaymath}
Thus the tangent space to $\mathcal{M}$ at $(E,\Phi)$ has a
decomposition
\begin{equation}\label{eq:tangent-dec}
  \HH^1(C^{\bullet}(E,\Phi))
    = \bigoplus_k \HH^1(C^{\bullet}_k(E,\Phi)).
\end{equation}

\begin{remark}\label{rem:duality-k}
  Using the definition of the $U_k$ and
  \eqref{eq:inner-invariant}, we have that
  \begin{displaymath}
    U^{\pm}_k \cong U^{\pm,*}_{-k}
  \end{displaymath}
  under the duality \eqref{eq:duality}.  Moreover, writing
  \begin{displaymath}
    \ad(\Phi)^{\pm}_k
      = \ad(\Phi)_{|U^{\pm}_k}\colon U^{\pm}_k \to
      U^{\mp}_{k+1}\otimes K,
  \end{displaymath}
  we have
  \begin{displaymath}
    \ad(\Phi)^{\pm}_{k,t} = (\ad(\Phi)^{\mp}_{-k-1})\otimes 1_{K^{-1}}.
  \end{displaymath}
\end{remark}

The calculations of Hitchin \cite[\S 8]{H2} show that eigenvalues
of the Hessian of $f$ at a critical point can be calculated as
follows.
\begin{proposition}\label{prop:negative-normal-bundle}
  Let $(E,\Phi)$ be a stable $G$-Higgs bundle which corresponds to a
  critical point of $f$, for $G=\GL(n,\C)$ or $G=\U(p,q)$.  In the
  decomposition \eqref{eq:tangent-dec} the eigenvalue $-k$ subspace
  for the Hessian of $f$ is isomorphic to
  $\HH^1(C^{\bullet}_k(E,\Phi))$.  In particular, the negative
  eigenspace at $(E,\Phi)$ for the Hessian is given by
  \begin{displaymath}
    \nu_{(E,\Phi)}^-(\mathcal{M}_l) \cong
      \bigoplus_{k>0} \HH^1(C^{\bullet}_k).
  \end{displaymath}
\end{proposition}

\begin{lemma}\label{lem:hyper-vanishing}
  Let $(E,\Phi)$ be a stable $G$-Higgs bundle which corresponds to a
  critical point of $f$.  Then
  \begin{displaymath}
    \HH^0(C^{\bullet}_k(E,\Phi)) = 0 \qquad\text{and}\qquad
    \HH^2(C^{\bullet}_k(E,\Phi)) = 0
  \end{displaymath}
  for $k>0$.
\end{lemma}

\begin{proof}
  {}From \eqref{eq:complex-dec} we have a decomposition
  \begin{displaymath}
    \HH^0(C^{\bullet}(E,\Phi))
      = \bigoplus_k \HH^0(C^{\bullet}_k(E,\Phi)).
  \end{displaymath}
  and we know from Proposition~\ref{prop:duality} that the only trivial
  endomorphisms of $(E,\Phi)$ are the scalars, which have weight zero
  in this decomposition.  This gives the vanishing of $\HH^0$.

  For the vanishing of $\HH^2$, consider first the case $G=\GL(n,\C)$.
  Then $U^+_k = U^-_k$ and, using Remark~\ref{rem:duality-k}, we see
  that the dual complex of $C^{\bullet}_k(E,\Phi)$ is isomorphic to
  the complex
  \begin{displaymath}
      C^{\bullet}_{-k-1}(E,\Phi) \otimes K^{-1} \colon
      U_{-k-1}^+\otimes K^{-1} \xrightarrow{-\ad(\Phi)} U_{-k}^- .
  \end{displaymath}
  The change in sign of $\ad(\Phi)$ does not influence the cohomology
  and hence Serre duality for hypercohomology gives
  \begin{displaymath}
    \HH^{2}(C^{\bullet}_{k}(E,\Phi))
      \cong \HH^{0}(C^{\bullet}_{-k-1}(E,\Phi))^*.
  \end{displaymath}
  It follows that $\HH^2(C^{\bullet}_{k}(E,\Phi))$ vanishes for $k \neq -1$.
  The case $G = \U(p,q)$ follows essentially from this, by using the
  fact that stability as a $\U(p,q)$-Higgs bundle implies stability as
  a $\GL(n,\C)$-Higgs bundle (see \cite[Proposition~3.19]{BGG1} for a
  detailed argument).
\end{proof}

\begin{proposition}\label{prop:morse-index}
  Let $(E,\Phi)$ be a stable $G$-Higgs bundle which corresponds to a
  critical point of $f$.  Then the Morse index of the corresponding
  critical submanifold $\mathcal{M}_l$ is
  \begin{displaymath}
    \Index(\mathcal{M}_l) =
      2 \sum_{k>0} \dim \HH^1(C^{\bullet}_{k}(E,\Phi)),
  \end{displaymath}
  where
  \begin{displaymath}
    \dim \HH^1(C^{\bullet}_{k}(E,\Phi))
    = -\chi(C^{\bullet}_{k}(E,\Phi)).
  \end{displaymath}
\end{proposition}

\begin{proof}
  This is immediate from Proposition~\ref{prop:negative-normal-bundle}
  and the vanishing of Lemma~\ref{lem:hyper-vanishing} (note that the
  Morse index is the real dimension of the $\sum \HH^1$, hence the
  factor of $2$).
\end{proof}

The following lemma is essentially Proposition~4.14 of \cite{BGG1}.
We provide a complete proof, taking this opportunity to correct some
inaccuracies in the argument given in \cite{BGG1}.

\begin{lemma}\label{lem:adjoint-criterion}
  Let $(E,\Phi)$ be a stable $G$-Higgs bundle which corresponds to a
  critical point of $f$, for $G=\GL(n,\C)$ or $G=\U(p,q)$.  Then
  \begin{displaymath}
     \chi(C^{\bullet}_k(E,\Phi)) \leq
     (g-1)\bigl(2 \rk(\ad(\Phi)_{k}^{+})
       - \rk(U_{k}^{+}) - \rk(U_{k+1}^{-}) \bigr).
  \end{displaymath}
  Furthermore, the vanishing $\chi(C^{\bullet}_k(E,\Phi)) = 0$ holds
  if and only if $\ad(\Phi)_k^+ \colon U_{k}^{+} \to U_{k+1}^{-}
  \otimes K$ is an isomorphism.
\end{lemma}

\begin{proof}
  In the following we shall use the abbreviated notations $C_k^{\bullet} =
  C^{\bullet}_k(E,\Phi)$ and
  \begin{displaymath}
    \Phi_k^{\pm} =
    \ad(\Phi)_k^{\pm}\colon U_k^+ \to U_{k-1}^-\otimes K.
  \end{displaymath}

  By the Riemann--Roch formula we have
  \begin{equation}\label{eq:RR-C_k}
    \chi(C_k^{\bullet}) = (1-g)\bigl(\rk(U_k^+)+\rk(U_{k+1}^-)\bigr)
      +\deg(U_k^+) - \deg(U_{k+1}^-),
  \end{equation}
  thus we can prove the inequality stated in the Lemma by
  estimating the difference $\deg(U_k^+) - \deg(U_k^-)$.  In order to
  do this, we note first that there are short exact sequences of sheaves
  \begin{displaymath}
    0 \to \ker(\Phi_k^+) \to U_k^+ \to \im(\Phi_k^+) \to 0
  \end{displaymath}
  and
  \begin{displaymath}
    0 \to \im(\Phi_k^+) \to U_{k+1}^-\otimes K
      \to \coker(\Phi_k^+)\to 0.
  \end{displaymath}
  It follows that
  \begin{equation}\label{eq:degUU1}
    \deg(U_k^+) - \deg(U_{k+1}^-)
    = \deg(\ker(\Phi_k^+)) + (2g-2)\rk(U_{k+1}^-) - \deg(\coker(\Phi_k^+)).
  \end{equation}
  We shall prove the following inequalities below.
  \begin{align}
    \deg(\ker(\Phi_k^+)) &\leq 0, \label{eq:3}\\
    - \deg(\coker(\Phi_k^+)) &\leq
    (2g-2)\bigl(-\rk(U_{k+1}^-) + \rk(\Phi_k^+)). \label{eq:4}
  \end{align}
  Combining these with \eqref{eq:degUU1} we obtain
  \begin{displaymath}
    \deg(U_k^+) - \deg(U_{k+1}^-)
    \leq (2g-2)\rk(\Phi_k^+),
  \end{displaymath}
  which, together with \eqref{eq:RR-C_k}, proves the inequality
  stated in the Lemma.

  It remains to prove \eqref{eq:3} and \eqref{eq:4}. For this we use
  the fact that the adjoint Higgs bundle $(U_G,\ad(\Phi))$ is
  semistable (one way of seeing this is to note that it supports a
  solution to Hitchin's equations).  Clearly, the subbundle
  $\ker(\Phi_k^+) \subseteq U_G$ is $\ad(\Phi)$-invariant and hence
  \begin{displaymath}
    \deg(\ker(\Phi_k^+)) \leq \deg(U_G) = 0,
  \end{displaymath}
  thus proving \eqref{eq:3}.

  In order to prove \eqref{eq:4} a bit more work needs to be done.
  Consider the dual of $\Phi_k^+$,
  \begin{displaymath}
    \Phi_k^{+,t}\colon U_{k+1}^{-,*}\otimes K^{-1} \to U_k^{+,*},
  \end{displaymath}
  and note that the image of $\Phi_k^+$ goes to zero under the
  restriction map
  \begin{displaymath}
    U_{k+1}^{-}\otimes K \to \ker(\Phi_k^{+,t})^*.
  \end{displaymath}
  Hence there is an induced map
  \begin{displaymath}
    \coker(\Phi_k^+) \to \ker(\Phi_k^{+,t})^*
  \end{displaymath}
  which is generically an isomorphism --- in fact, its kernel is the
  torsion subsheaf of $\coker(\Phi_k^+)$. It follows that
  \begin{displaymath}
    \deg(\coker(\Phi_k^+)) \geq \deg(\ker(\Phi_k^{+,t})^*).
  \end{displaymath}
  Since $\ker(\Phi_k^{+,t})$ is locally free (in fact a subbundle)
  this shows that
  \begin{equation}\label{eq:1}
    -\deg(\coker(\Phi_k^+)) \leq \deg(\ker(\Phi_k^{+,t})),
  \end{equation}
  the difference being the degree of the torsion subsheaf of
  $\coker(\Phi_k^+)$.  Now Remark~\ref{rem:duality-k} tells us that we
  have a commutative diagram
  \begin{displaymath}
    \begin{CD}
      U_{k+1}^{-,*}\otimes K^{-1} @>\Phi_k^{+,t}>> U_k^{+,*} \\
      @V{\cong}VV  @V{\cong}VV \\
      U_{-k-1}^-\otimes K^{-1} @>{-\Phi_{-k-1}^-\otimes 1_{K^{-1}}}>>
        U_{k}^+, \\
    \end{CD}
  \end{displaymath}
  and thus
  \begin{displaymath}
    \ker(\Phi_k^{+,t}) \cong \ker(\Phi_{-k-1}^-)\otimes K^{-1}
  \end{displaymath}
  from which we conclude that
  \begin{displaymath}
    \deg(\ker(\Phi_k^{+,t})) = \deg(\ker(\Phi_{-k-1}^-))
      - (2g-2)\rk(\ker(\Phi_{-k-1}^-)).
  \end{displaymath}
  Again we apply semistability of $(U_G,\ad(\Phi))$ to the
  $\ad(\Phi)$-invariant subbundle $\ker(\Phi_{-k-1}^-)$ to obtain
  \begin{equation}\label{eq:d-ker-phi-t}
    \deg(\ker(\Phi_k^{+,t})) \leq -(2g-2)\rk(\ker(\Phi_{-k-1}^-)).
  \end{equation}
  But clearly, $\rk(\Phi_k^+) = \rk(\Phi_k^{+,t}) = \rk(\Phi_{-k-1}^-)$
  and $\rk(U_{k+1}^-) = \rk(U_{-k-1}^{-,*}) = \rk(U_{-k-1}^{-})$
  so
  $$\rk(\ker(\Phi_{-k-1}^-)) = \rk(U_{k+1}^-) - \rk(\Phi_k^+).$$
  Combining this fact with \eqref{eq:1} and
  \eqref{eq:d-ker-phi-t} concludes the proof of \eqref{eq:4}.

  Finally, to prove the last statement of the Lemma we note that if
  $\chi(C_k^{\bullet}) = 0$ then $\rk(\Phi_k^+) = \rk(U_k^+) =
  \rk(U_{k+1}^- \otimes K)$ and equality holds in \eqref{eq:1}, thus
  showing that $\Phi_k^+$ is an isomorphism.
\end{proof}

\begin{proposition}\label{index-bound}
  \begin{enumerate}
  \item\label{item:index-upq}
  For ${\mathcal M}={\mathcal
    M}(p,q,a,b)$
$$
\Index(\mathcal{M}_l)\geq 2g-2
$$
for every non-minimal critical submanifold ${\mathcal M}_l\subset
\mathcal{M}$.
  \item\label{item:index-gln}
  For ${\mathcal M}=\mathcal{M}(n,d)$
$$
\Index(\mathcal{M}_l)\geq (n-1)(2g-2)
$$
for every non-minimal critical submanifold ${\mathcal M}_l\subset
\mathcal{M}$.
  \end{enumerate}
\end{proposition}

\begin{proof}
  (\ref{item:index-upq}) Let $k_0$ be the largest $k$ such that
  $\chi(C_k^\bullet(E,\Phi)) \neq 0$. Since $\mathcal{M}_l$ is
  non-minimal, by Proposition~\ref{prop:morse-index} we have $k_0 >
  0$. The proof of \cite[Proposition~4.17]{BGG1} shows that the
  restriction of $\ad(\Phi)_k^+ \colon U_{k}^{+} \to U_{k+1}^{-}
  \otimes K$ to a fibre is never an isomorphism (in the notation of
  that proof, $k_0=m-1$).  Hence the right hand side of the inequality
  of Lemma~\ref{lem:hyper-vanishing} is strictly negative.  Now the
  result follows from this inequality and
  Proposition~\ref{prop:morse-index}.

  (\ref{item:index-gln}) We recall (cf.\ \cite{H,H2,S2}) that the
  decomposition $U = \bigoplus U_k$ comes from a decomposition $E =
  E_1 \oplus \dots \oplus E_m$ with $U_k =
  \bigoplus_{k=j-i}\Hom(E_i,E_j)$.  In particular, the weights $k$ are
  consecutive integers.  Thus Proposition~\ref{prop:morse-index},
  together with Riemann--Roch and the fact that $U^+_{G} = U_{G}^-$
  for $G=\GL(n,\C)$, gives
  \begin{align*}
    \tfrac{1}{2}\Index(\mathcal{M}_l) &=
      (g-1)\sum_{k\geq 1} \bigl(\rk(U_k) + \rk(U_{k+1}) -
        2\rk(\ad(\Phi)_k\bigr) \\
      &= (g-1)\bigl(\rk(U_1) + 2\rk(U_{k \geq 2}) - 2\rk(\ad(\Phi)_{k
        \geq 1})\bigr).
  \end{align*}
  But clearly the rank of $\ad(\Phi)_{k \geq 1}\colon U_{k \geq 1} \to
  U_{k \geq 2}\otimes K$ is less than or equal to the rank of $U_{k
    \geq 2}$ and hence
  \begin{displaymath}
    \tfrac{1}{2}\Index(\mathcal{M}_l) \geq (g-1)\rk(U_1).
  \end{displaymath}
  Let $\nu_i = \rk(E_i)$.  Then $\sum \nu_i =n$ and $\rk(U_1) =
  \nu_1\nu_2 + \dots + \nu_{m-1}\nu_m$.  One easily shows that
  $\nu_1\nu_2 + \dots + \nu_{m-1}\nu_m \geq n-1$.  This finishes the
  proof of (\ref{item:index-gln}).
\end{proof}

\begin{remark}\label{rem:referee}
  Our initial estimate in (\ref{item:index-gln}) of
  Proposition~\ref{index-bound} was $\Index(\mathcal{M}_l)\geq 2g-2$.
  It was pointed out to us by an anonymous referee that this could be
  improved, and also that an alternative way of proving this estimate
  is as follows.  The absolute minimum of $f$ on $\mathcal{M}(n,d)$ is
  $M(n,d)$, so $H^*(M(n,d))$ injects into $H^*(\mathcal{M}(n,d))$
  because $f$ is perfect. For the same reason, any critical
  submanifold of index $l$ gives a non-trivial contribution to the
  cohomology of $\mathcal{M}(n,d)$ in dimension $l$, which is not in
  the image of $H^*(M(n,d))$.  Now Markman \cite{markman:2002} shows
  that $H^*(B\bar{\mathcal{G}})$ (the cohomology of the classifying
  space of the reduced gauge group) surjects onto
  $H^*(\mathcal{M}(n,d))$. On the other hand, Uhlenbeck--Daskalopoulos
  \cite{DU} prove that
  $H^r(B\bar{\mathcal{G}})$ is isomorphic to $H^r(M(n,d))$ for
  $r<(2g-2)(n-1)$.  Hence no critical submanifold can have index $l <
  (2g-2)(n-1)$.
\end{remark}

\subsection{Local Minima}
\label{sec:local-minima}

The minima of the Morse function on ${\mathcal M}(n,d)$ is given by
the following \cite{H}.

\begin{proposition} \label{minima-gl}
Let $ \mathcal{N}(n,d)\subset {\mathcal M}(n,d)$ be the set of local minima.
Then
$$
 \mathcal{N}(n,d)=\{ (E,\Phi) \in \mathcal{M}(n,d)\suchthat \Phi=0\}.
$$
Hence  $\mathcal{N}(n,d)$ coincides with $M(n,d)$, the moduli space of
semistable vector bundles of rank $n$ and degree $d$, which equals
the subvariety $M^s(n,d)\subset M(n,d)$ corresponding to stable bundles if
$\GCD(n,d)=1$.

\end{proposition}

The minima of the Morse function on ${\mathcal M}(p,q,a,b)$ have been
characterized in \cite{BGG1}. One has the following.

\begin{proposition} \label{minima-upq}
Let $ \mathcal{N}(p,q,a,b)\subset {\mathcal M}(p,q,a,b)$ be the set of local minima.
Then
$$
 \mathcal{N}(p,q,a,b)=\{ (E,\Phi) \in \mathcal{M}(p,q,a,b)\suchthat
 \beta= 0 \;\;\mbox{or}\;\;\gamma=0\}.
$$
More precisely, let $(E,\Phi)\in \mathcal{N}(p,q,a,b)$.  Then
  \begin{enumerate}

  \item[$(1)$] $\beta = 0$ if and only if $a/p > b/q$
  (i.e.\ $\tau > 0$).
  \item[$(2)$] $\gamma = 0$ if and only if $a/p < b/q$
(i.e.\ $\tau < 0$).
\end{enumerate}
\end{proposition}

\begin{remark}
Since we are assuming $\GCD(p+q,a+b)=1$ then $\tau\neq 0$.
\end{remark}

\section{Homotopy groups}

\subsection{Homotopy groups of $\mathcal{M}(n,d)$}

Combining Propositions \ref{homotopy-minima}, \ref{index-bound}
and \ref{minima-gl} we have the following.

\begin{theorem}\label{homotopy-minima-gl} Let $\GCD(n,d)=1$. Then
$$
\pi_i(\mathcal{M}(n,d))\cong \pi_i(M(n,d)),
\;\;\;{for}\;\;\; i\leq 2(g-1)(n-1)-2.
$$
\end{theorem}
Now, the homotopy groups of $M(n,d)$ have been computed
by  Daskalopoulos and Uhlenbeck \cite{DU} (here  $n$ and $d$ are not assumed
to be coprime). Their result is the following.

\begin{theorem}\label{daskalopoulos-uhlenbeck} Let $M^s(n,d)$ be the
moduli space of stable vector bundles of rank $n$ and degree $d$. Assume
that $n>1$ and $(n,g)\neq (2,2)$. Then
\begin{enumerate}
\item $\pi_1(M^s(n,d))\cong H_1(X,\Z)$;
\item $\pi_2(M^s(n,d))\cong \Z \oplus \Z_\GCD(n,d)$;
\item $\pi_i(M^s(n,d))\cong \pi_{i-1}(\mathcal{G}), \;\;\;
\mbox{for}\;\;\; 2<i\leq 2(g-1)(n-1)-2$, where $\mathcal{G}$ is the
unitary gauge group.
\end{enumerate}
\end{theorem}

\begin{remark}
  The proof of (1) when $n$ and $d$ are coprime is given by
  Atiyah--Bott \cite{AB}.
\end{remark}

As a corollary of Theorems \ref{homotopy-minima-gl} and
\ref{daskalopoulos-uhlenbeck} we have the following.

\begin{theorem}\label{homotopy-higgs-n-d}
Assume that $n>1$ and $\GCD(n,d)=1$ and let $g\geq 3$. Then
\begin{enumerate}
\item $\pi_1(\mathcal{M}(n,d))\cong H_1(X,\Z)$;
\item $\pi_2(\mathcal{M}(n,d))\cong \Z$;
\item\label{item:3}
  $\pi_i(\mathcal{M}(n,d))\cong \pi_{i-1}(\mathcal{G}), \;\;\;
  \mbox{for}\;\;\; 2<i\leq 2(g-1)(n-1)-2$.
\end{enumerate}
\end{theorem}

\begin{remark}
  As a consequence of Theorem \ref{homotopy-minima-gl} and the
  connectedness of $M(n,d)$ \cite{NS} one obtains that
  $\mathcal{M}(n,d)$ is also connected \cite{H,S3}.

A proof of (1) when $n=2$ is given by Hitchin \cite{H}.
\end{remark}

\begin{remark}
  When $n=2$ Hausel \cite[Theorem 7.5.7]{hausel:1998a} proved that the
  isomorphism (\ref{item:3}) holds for $i \leq 4g-8$, which is twice
  as good as our estimate.  It would be very interesting to see if
  this result can be generalized to higher $n$.
\end{remark}

\subsection{Moduli space of triples}

The next step is to identify the spaces $\mathcal{N}(p,q,a,b)$ as
moduli spaces in their own right. They turn out to be examples of the
moduli spaces of triples studied in \cite{BG} \cite{BGG1} and
\cite{BGG2}.  We briefly recall the relevant definitions and results.
See \cite{BGG1} for details.

A \emph{holomorphic triple} on $X$, $T = (E_{1},E_{2},\phi)$, consists
of two holomorphic vector bundles $E_{1}$ and $E_{2}$ on $X$ and a
holomorphic map $\phi \colon E_{2} \to E_{1}$.  Denoting the ranks
$E_1$ and $E_2$ by $n_1$ and $n_2$, and their degrees by $d_1$ and
$d_2$, we refer to $(n_1,n_2,d_1,d_2)$ as the {\em type} of the
triple.

A homomorphism from $T' = (E_1',E_2',\phi')$ to $T = (E_1,E_2,\phi)$
is a commutative diagram
\begin{displaymath}
  \begin{CD}
    E_2' @>\phi'>> E_1' \\
    @VVV @VVV  \\
    E_2 @>\phi>> E_1.
  \end{CD}
\end{displaymath}
$T'=(E_1',E_2',\phi')$ is a subtriple of $T = (E_1,E_2,\phi)$
if the homomorphisms of sheaves $E_1'\to E_1$ and $E_2'\to E_2$ are
injective.

For any $\alpha \in \R$ the \emph{$\alpha$-degree} and
\emph{$\alpha$-slope} of $T$ are
defined to be
\begin{align*}
  \deg_{\alpha}(T)
  &= \deg(E_{1}) + \deg(E_{2}) + \alpha
  \rk(E_{2}), \\
  \mu_{\alpha}(T)
  &=
  \frac{\deg_{\alpha}(T)}
  {\rk(E_{1})+\rk(E_{2})} \\
  &= \mu(E_{1} \oplus E_{2}) +
  \alpha\frac{\rk(E_{2})}{\rk(E_{1})+
    \rk(E_{2})}.
\end{align*}
The triple $T = (E_{1},E_{2},\phi)$ is
\emph{$\alpha$-stable} if
\begin{equation}\label{eqn:alpha-stability}
  \mu_{\alpha}(T')
  < \mu_{\alpha}(T)\
\end{equation}
for any proper sub-triple $T' = (E_{1}',E_{2}',\phi')$. Define
\emph{$\alpha$-semistability} by replacing (\ref{eqn:alpha-stability})
with a weak inequality. A triple is called
\emph{$\alpha$-polystable} if it is the direct sum of $\alpha$-stable
triples of the same $\alpha$-slope. It is {\it strictly
$\alpha$-semistable (polystable)} if it is $\alpha$-semistable (polystable)
but not $\alpha$-stable.

We denote the moduli space of isomorphism classes of $\alpha$-{\it
polystable} triples of type $(n_1,n_2,d_1,d_2)$ by
\begin{equation}\label{eqn:N(n,d)}
  \mathcal{N}_\alpha
  = \mathcal{N}_\alpha(n_1,n_2,d_1,d_2).
\end{equation}
Using Jordan--H\"older filtrations of $\alpha$-semistable triples
one can  define $S$-equivalence, and view
$  \mathcal{N}_\alpha$ as  the moduli space of $S$-equivalence classes of
$\alpha$-semistable triples. The isomorphism
classes of $\alpha$-{\it stable} triples form a subspace which we
denoted by $\mathcal{N}_\alpha^s$.

\begin{proposition}[\cite{BG}]
  \label{prop:alpha-range}
The moduli space   $\mathcal{N}_\alpha(n_1,n_2,d_1,d_2)$ is a complex
projective variety. A
necessary condition for the moduli space
$\mathcal{N}_\alpha(n_1,n_2,d_1,d_2)$ to be non-empty is
\begin{equation}
\begin{cases}
0\leq \alpha_m \leq \alpha \leq \alpha_M &\text{if $n_1\neq n_2$}\\
0\leq \alpha_m \leq \alpha &\text{if $n_1= n_2$}
\end{cases}
\end{equation}
where
\begin{align}
 \alpha_m &= \mu_1-\mu_2, \label{alpha-bounds-m} \\
      \alpha_M &=
      (1+ \frac{n_1+n_2}{|n_1 - n_2|})(\mu_1 - \mu_2)
  \label{alpha-bounds-M}
\end{align}
and $\mu_1=\frac{d_1}{n_1}$, $\mu_2=\frac{d_2}{n_2}$.
\end{proposition}

Whenever necessary we shall
indicate the dependence of $\alpha_m$ and $\alpha_M$ on $(n_1,n_2,d_1,d_2)$
by writing $\alpha_m=\alpha_m(n_1,n_2,d_1,d_2)$, and similarly for $\alpha_M$.

Within the allowed range for $\alpha$ there is a discrete set of {\it
critical values}. These are the values of
$\alpha$ for which it is numerically possible to have a subtriple
$T'=(E_1',E_2',\phi')$ such that $\mu(E'_1 \oplus E'_2)\ne\mu(E_1
\oplus E_2)$ but
$\mu_{\alpha}(T')=\mu_{\alpha}(T')$. All other values of $\alpha$ are
called {\it generic}. The critical values of $\alpha$ are precisely
the values for $\alpha$ at which the stability properties of a triple
can change, i.e.\ there can be triples which are strictly
$\alpha$-semistable, but either
$\alpha'$-stable or $\alpha'$-unstable for $\alpha'\ne\alpha$.

The following result relates the stability conditions for holomorphic
triples and that for $\U(p,q)$-Higgs bundles.

\begin{proposition}
  \label{prop:triple-higgs-stability}
  A $\U(p,q)$-Higgs bundle $(E,\Phi)$ with $\beta =0$ or $\gamma=0$ is
  \mbox{(semi)}stable if and only if the corresponding holomorphic
  triple  is
  $\alpha$-(semi)stable for $\alpha = 2g-2$.
\end{proposition}

Combining Propositions \ref{minima-upq} and
\ref{prop:triple-higgs-stability}, we have the following
characterization of the subspace of local minima
$\mathcal{N}(p,q,a,b)$.
\begin{theorem}
  \label{thm:minima=triple-moduli}
  Let $\mathcal{N}(p,q,a,b)$ be the subspace of local minima of
  $f$ on $\mathcal{M}(p,q,a,b)$ and let $\tau$ be the Toledo
  invariant.

  If $a/p < b/q$, or equivalently if $\tau < 0$,
   then $\mathcal{N}(p,q,a,b)$ can be
  identified with the moduli space of $\alpha$-polystable triples of type
  $(p,q, a + p(2g-2), b)$, with $\alpha = 2g-2$.

  If $a/p > b/q$,  or equivalently if $\tau > 0$,
  then $\mathcal{N}(p,q,a,b)$ can be
  identified with the moduli space of $\alpha$-polystable triples of type
  $(q,p, b + q(2g-2),a)$, with $\alpha = 2g-2$.

  That is,
  \begin{displaymath}
    \mathcal{N}(p,q,a,b)
    \cong \begin{cases}
    \mathcal{N}_{2g-2}(p,q,a + p(2g-2),b)
    &\text{if $a/p <b/q$ (equivalently $\tau < 0$)}\\
    \mathcal{N}_{2g-2}(q,p,b + q(2g-2),a)
    &\text{if $a/p > b/q$ (equivalently $\tau >0$)}\end{cases}
  \end{displaymath}
\end{theorem}

In view of Theorem~\ref{thm:minima=triple-moduli}
it is important to
understand where $2g-2$ lies in relation to the range (given by
Proposition \ref{prop:alpha-range}) for the stability parameter
$\alpha$. One has the following.

\begin{proposition}\label{prop:summary-range-2g-2} Fix $(p,q,a,b)$.
Then
\begin{equation}
0\le|\tau|\le\min\{p,q\}(2g-2)
\Leftrightarrow
       0<\alpha_m(p,q,a,b)\le 2g-2 \le
       \alpha_M(p,q,a,b) \text{if $p\ne q$}
\end{equation}
\end{proposition}

Proposition \ref{prop:summary-range-2g-2} shows that in order to study
$\mathcal{N}(p,q,a,b)$ for different values of the Toledo invariant, we
need to understand the moduli spaces of triples for values of $\alpha$
that may lie anywhere (including at the extremes $\alpha_m$ and
$\alpha_M$) in the $\alpha$-range given in
Proposition~\ref{prop:alpha-range}.

We can assume   $n_1>n_2$, since by  triples
duality one has the following.

\begin{proposition} \label{triples-duality}
$\mathcal{N}_\alpha(n_1,n_2,d_1,d_2)\cong
\mathcal{N}_\alpha(n_2,n_1,-d_2,-d_1)$.
\end{proposition}

Recall that the allowed range for the stability parameter is
$\alpha_m\leq\alpha\leq\alpha_M$,
where $\alpha_m=\mu_1-\mu_2$ and
$\alpha_M=\frac{2n_1}{n_1-n_2}\alpha_m$, and we assume  that
$\mu_1-\mu_2>0$.   We describe the moduli space $\mathcal{N}_\alpha$
for $2g-2 \leq \alpha < \alpha_M$.

Let $\alpha_L$ be the largest critical value in
 $(\alpha_m,\alpha_M)$, and let $\mathcal{N}_L$
(respectively $\mathcal{N}^s_L$) denote the moduli space of
$\alpha$-polystable (respectively
$\alpha$-stable) triples for
$\alpha_L<\alpha<\alpha_M$. We refer to $\mathcal{N}_L$ as the `large
$\alpha$' moduli space.

\begin{proposition}
\label{prop:triple-stable-implies-bundles-semistable}
Let $T=(E_1,E_2,\phi)$ be an  $\alpha$-semistable triple for
some $\alpha$\ in the range $\alpha_L<\alpha<\alpha_M$. Then
$T$ is of the form
\begin{equation*}
0 \lto E_2 \overset{\phi}{\lto} E_1 \lto F \lto 0,
\end{equation*}
with  $F$ locally free, and $E_2$ and $F$ are semistable.

\end{proposition}
In the converse direction we have\footnote{This Proposition replaces
  Proposition~7.6 of \cite{BGG2}. We thank Stefano Pasotti, Francesco
  Prantil and Carlos Tejero for pointing out errors in this Proposition.}:
\begin{proposition}
Let $T=(E_1,E_2,\phi)$ be a triple of  the form
\begin{equation}\label{extension} 0 \lto E_2 \overset{\phi}{\lto} E_1 \lto F
\lto 0,
\end{equation}
with $F$ locally free and such that the extension is non-trivial.  If $E_2$
 and $F$ are stable then $T$ is $\alpha$-stable for $\alpha$ in the
range $\alpha_L<\alpha<\alpha_M$.
\end{proposition}

\begin{proof}
  Regarding the top line of the diagram
  \begin{equation} \label{extension2}
     \begin{CD}
         E_2 @>>> \phi(E_2) \\
         @VVV   @VVV  \\
         E_2 @>\phi>> E_1 \\
         @VVV   @VVV  \\
         0 @>>> F
    \end{CD}
  \end{equation}
  as a subtriple $T'$ of $T$, and the bottom line as a quotient triple
  $T''$, we can consider $T$ as an extension of triples
  \begin{displaymath}
    0 \lto T' \lto T \lto T'' \lto 0\ .
  \end{displaymath}
  It follows from stability of $E_2$ that the subtriple $T'$ is
  $\alpha$-stable for any $\alpha>0$. It follows similarly from
  stability of $F$, that the quotient triple $T''$ is $\alpha$-stable
  for any $\alpha$.  In particular $T'$ and $T''$ are both
  $\alpha_M$-stable.

  A simple calculation shows that
  \begin{equation}\label{eq:7}
    \mu_{\alpha_M}(T') = \mu(E_2) +\frac{1}{2}\alpha_M
    = \mu(F) = \mu_{\alpha_M}(T'')\ .
  \end{equation}
  It is a general fact that an extension of $\alpha$-semistable
  triples of the same $\alpha$-slope is itself $\alpha$-semistable.
  Thus we deduce from (\ref{eq:7}) that the triple $T$ is
  $\alpha_M$-semistable.

  It remains to show that $T$ is $\alpha$-stable for
  $\alpha_L<\alpha<\alpha_M$.  We do this by showing that there is no
  $\alpha$-destabilizing subtriple, i.e., a subtriple $S$ of $T$ such
  that $\mu_\alpha(S) \geq \mu_\alpha(T)$.

  To do this, we first observe that the extension (\ref{extension2})
  is a Jordan--H\"older filtration of $T$ considered as an
  $\alpha_M$-semistable object.  This follows since $T'$ and $T''$ are
  $\alpha_M$-stable and have the same $\alpha_M$-slope. Hence the
  associated graded object for $T$ in the category of
  $\alpha_M$-semistable triples is
  \begin{equation}\label{eq:2}
    \Gr(T) = T' \oplus T''\ .
  \end{equation}

  Now assume that $S \subseteq T$ is $\alpha$-destabilizing for
  $\alpha$ in the range $\alpha_L<\alpha<\alpha_M$.  By continuity of
  $\mu_\alpha(S)$ in $\alpha$, we have $\mu_{\alpha_M}(S) \geq
  \mu_{\alpha_M}(T)$ and hence $\alpha_M$-semistability of $T$ implies
  that $\mu_{\alpha_M}(S) = \mu_{\alpha_M}(T)$.  It follows that in
  $\Gr(T) =T'\oplus T''$  the triple $S$ must
  be isomorphic to either $T'$ or $T''$.

  We first show that $S$ cannot be isomorphic to $T'$, i.e, that the
  subtriple $T'$ is not destabilizing for $\alpha<\alpha_M$. The key
  piece of information is that
  \begin{align*}
    \frac{n_2(T')}{n(T')} &= \frac{n_2}{2n_2} = \frac{1}{2}\ , \\
    \frac{n_2(T)}{n(T)} &= \frac{n_2}{n_1+n_2} < \frac{1}{2}\ , \\
  \end{align*}
  where, for any triple $T = (E_2,E_1,\phi)$, we write $n_i(T) =
  \rk(E_i)$ and $n(T) = \rk(E_2 \oplus E_1)$.
  Hence
  \begin{equation}\label{eq:5}
    \frac{n_2(T')}{n(T')} > \frac{n_2(T)}{n(T)}>\frac{n_2(T'')}{n(T'')} \ .
  \end{equation}
  But, since $\mu_{\alpha_M}(T') = \mu_{\alpha_M}(T)$,
  \begin{displaymath}
    \mu_{\alpha}(T') - \mu_{\alpha}(T) =
    (\alpha-\alpha_M)\left(\frac{n_2(T')}{n(T')} -
      \frac{n_2}{n}\right)
    < 0
  \end{displaymath}
  for $\alpha < \alpha_M$.

  Finally we show that $S$ cannot be isomorphic to $T''$. In fact, if
  $T$ has a subtriple isomorphic to $T''$, then $E_1$ has a subbundle,
  $\tilde{F}$, isomorphic to $F$. The composition of the isomorphism
  from $F$ to $\tilde{F}$ with the projection from $E_1$ to $F$
  produces a homomorphism
  $$\psi:F\to F \ .$$
  Since $F$ is stable, $\psi$ is either zero or a multiple of the
  identity. If it is zero, then there must be a non-trivial
  homomorphism from $F$ to $E_2$.  This is impossible since
  $\mu(\tilde{F})>\mu(E_2)$, and both are stable bundles. Hence the
  the isomorphism from $F$ to $\tilde{F}$ splits the extension
  (\ref{extension}). But this contradicts our assumptions.
\end{proof}

\begin{theorem}\label{thm:largealpha}
Assume that $n_1>n_2$ and $d_1/n_1> d_2/n_2$. Then the moduli space
$\mathcal{N}^s_L=\mathcal{N}^s_L(n_1,n_2,d_1,d_2)$
is  smooth
of dimension
$$
(g-1)(n_1^2 + n_2^2 - n_1 n_2) - n_1 d_2 + n_2 d_1 + 1,
$$
and includes  a $\mathbb{P}^N$-fibration
$\mathcal{P}$
 over
$M^s(n_1-n_2,d_1-d_2) \times M^s(n_2,d_2)$, where $M^s(n,d)$ is the moduli
space of stable bundles of rank $n$ and degree $d$, and
$N=n_2d_1-n_1d_2+n_1(n_1-n_2)(g-1)-1$. Moreover, the complex codimension of
$\mathcal{N}^s_L\setminus \mathcal{P}$ is at least $g-1$. In particular,
$\mathcal{N}_L^s(n_1,n_2,d_1,d_2)$ is non-empty and irreducible.

If $\GCD(n_1-n_2,d_1-d_2)=1$ and $\GCD(n_2,d_2)=1$, then
$\mathcal{N}_L^s(n_1,n_2,d_1,d_2)$ is isomorphic to $\mathcal{P}$.
\end{theorem}
\begin{proof}
The birational equivalence between $\mathcal{P}$ and  $\mathcal{N}^s_L$
is proved in \cite{BGG2}. To obtain the precise estimate of the
codimension of $\mathcal{N}^s_L\setminus \mathcal{P}$
in $\mathcal{N}^s_L$ we see that,
by Proposition \ref{prop:triple-stable-implies-bundles-semistable},
it suffices to
estimate the dimension of stable  triples like (\ref{extension})
with $E_2$ and $F$ semistable.

Now, for any family of semistable bundles
the complex  codimension of the  set of strictly semistable bundles
is at least $g-1$. A computation of the precise estimate can be found in
\cite{BGMMN}.
The proof is finished by observing that for a stable triple
of the form  (\ref{extension}) $H^0(X,E_2\otimes F^*)=0$
(see \cite{BGG2}).
\end{proof}

The following is proved in \cite{BGG2}.

\begin{theorem}\label{thm:irreducibility-moduli-stable-triples}
Let $\alpha$ be any value in the range $\alpha_m<2g-2\leq\alpha<
\alpha_M$. Then
$\mathcal{N}^s_\alpha$ is birationally equivalent to
$\mathcal{N}^s_L$. Moreover, they are isomorphic outside of a set of
complex codimension greater or equal than $g-1$.
In particular, $\mathcal{N}^s_\alpha$
is non-empty and irreducible.
\end{theorem}

\subsection{Homotopy groups of moduli spaces of triples}

The strategy to compute the homotopy groups of $\mathcal N(p,q,a,b)$
is to compute first those of the moduli space of $\alpha$-stable
triples $\mathcal{N}^s_\alpha$ for large $\alpha$.

Let $n_1>n_2$ and let $\mathcal{P}\subset \mathcal{N}_L$ be the
$\PP^N$-fibration over $M^s(n_1-n_2,d_1-d_2)\times   M^s(n_2,d_2)$ given in
Theorem \ref{thm:largealpha}. As a consequence of Theorems
\ref{thm:largealpha}  and \ref{thm:irreducibility-moduli-stable-triples} we
have the following.

\begin{proposition}\label{homotopy-large-alpha}
Let $2g-2\leq \alpha <\alpha_M$. Then
$$
\pi_i (\mathcal{N}^s_\alpha(n_1,n_2,d_1,d_2)) \cong \pi_i
(\mathcal{N}^s_L(n_1,n_2,d_1,d_2)) \cong \pi_i(\mathcal{P})\;\;\;
\mbox{for} \;\;\; i\leq 2g-4.
$$
\end{proposition}

Associated to the $\PP^N$-fibration $\mathcal{P}$ over
$M^s(n_1-n_2,d_1-d_2)\times M^s(n_2,d_2)$ there is a homotopy sequence

\begin{multline}\label{homotopy-sequence}
  \cdots  \lra  \pi_i(\PP^N) \lra  \pi_i(\mathcal{P})
  \lra \pi_i(M^s(n_2,d_2))\times  \pi_i(M^s(n_1-n_2,d_1-d_2)) \\
  \lra
\pi_{i-1}(\PP^N) \lra \cdots
\end{multline}

\begin{proposition}\label{homotopy-triples}
Let $n_1>n_2$ and $n_2d_1>n_1d_2$. Assume that $(n_2,g)\neq (2,2)$ and
$(n_1- n_2,g)\neq (2,2)$ (for our applications we will actually assume
$g\neq 3$). Then
\begin{enumerate}

\item $\pi_1(\mathcal{P})\cong \pi_1(M^s(n_2,d_2))\times
\pi_1(M^s(n_1-n_2,d_1-d_2))
\cong H_1(X,\Z) \oplus H_1(X,\Z)$;

\item $\pi_2(\mathcal{P})$ is the middle term in an exact sequence
\begin{equation*}
0\lra\Z\lra\pi_2(\mathcal{P})\lra Q(n_1,n_2,d_1,d_2)\lra 0\
\end{equation*}

\noindent where

\begin{equation*}
Q(n_1,n_2,d_1,d_2)=
\begin{cases} \Z \oplus \Z \oplus \Z_{\GCD(n_2,d_2)}
\oplus\Z_{\GCD(n_1-n_2,d_1-d_2)}&\text{if $n_2>1$\ and $(n_1-n_2)>1$}\\
\Z \oplus \Z_{\GCD(n_2,d_2)}&\text{if $n_2>1$\ and $(n_1-n_2)=1$}\\
\Z\oplus\Z_{\GCD(n_1-n_2,d_1-d_2)}&\text{if $n_2=1$\ and $(n_1-n_2)>1$}\\
0&\text{if $n_2=1$\ and $n_1=2$}\\
\end{cases}
\end{equation*}

\end{enumerate}
\end{proposition}

\begin{remark}
  \label{rem:coprime-extension}
  It follows immediately from the exact sequence in (2) of
  Proposition~\ref{homotopy-triples} that the free part of the
  finitely generated abelian group $\pi_2(\mathcal{P})$ equals the
  direct sum of $\Z$ and the free part of $Q(n_1,n_2,d_1,d_2)$.  In
  particular, if the co-primality conditions $\GCD(n_2,d_2) = 1$ and
  $\GCD(n_1-n_2,n_2-d_2) = 1$ hold, then we have a complete
  description of $\pi_2(\mathcal{P})$ as the direct sum $\Z \oplus
  Q(n_1,n_2,d_1,d_2)$.  Also, under any circumstances, it follows that
  $\pi_2(\mathcal{P}) \otimes \Q \cong \Q \oplus Q(n_1,n_2,d_1,d_2)
  \otimes \Q$ so, for rational homotopy, our results are complete.
\end{remark}

\begin{proof}[Proof of Proposition~\ref{homotopy-triples}]

{}From the  homotopy sequence (\ref{homotopy-sequence}), since
$\pi_0(\PP^N) =\pi_1(\PP^N)=0$,  we deduce that
$\pi_1(\mathcal{P})\cong \pi_1(M^s(n_2,d_2))\times
\pi_1(M^s(n_1-n_2,d_1-d_2))$.
Statement (1) follows from Theorem \ref{daskalopoulos-uhlenbeck}.

Since $\pi_1(\PP^N)=0$, (\ref{homotopy-sequence})  gives

\begin{equation}\label{2-homotopy}
 ...  \lra  \pi_2(\PP^N) \lra  \pi_2(\mathcal{P})  \lra
\pi_2(M^s(n_2,d_2))\times  \pi_2(M^s(n_1-n_2,d_1-d_2)) \lra 0.
\end{equation}

On the other hand, by Hurewicz' theorem
$\pi_2(\PP^N))\cong H_2(\PP^N,\Z)\cong \Z$. Now, the map
$f: \Z\cong \pi_2(\PP^N)   \lra \pi_2(\mathcal{P})$  in (\ref{2-homotopy})
is injective since
one has the commutative diagram
\begin{displaymath}
  \begin{CD}\Z\cong
    \pi_2(\PP^N)     @>f>> \pi_2(\mathcal{P})\\
       @V\|VV @VVV  \\
    H_2(\PP^N,\Z) @>>> H_2(\mathcal{P},\Z).
  \end{CD}
\end{displaymath}
and $H_2(\PP^N,\Z)\lra H_2(\mathcal{P},\Z)$ must be injective because the
restriction of an ample line bundle over $\mathcal{P}\subset \mathcal{N}_L$
to $\PP^N$ must give an ample line bundle.  Note that the natural map
$H_2(\PP^N)\lra H_2(\mathcal{N}_L)$ is injective and factors through
$H_2(\PP^N)\lra H_2(\mathcal{P})\lra H_2(\mathcal{N}_L)$.  Now, we obtain
(2) from Theorem \ref{daskalopoulos-uhlenbeck} and the fact that if $n=1$
then the moduli space $M^s(n,d)$ is the Jacobian of degree $d$ line
bundles.
\end{proof}

As a corollary of Proposition \ref{homotopy-large-alpha}
and Proposition \ref{homotopy-triples} we have the following.

\begin{theorem}\label{thm:homotopy-triples}
Assume $n_1>n_2$, $n_2d_1>n_1d_2$, $g\geq 3$, and $2g-2\leq \alpha
<\alpha_M$. Then

\begin{enumerate}

\item $\pi_1(\mathcal{N}^s_\alpha(n_1,n_2,d_1,d_2)) \cong H_1(X,\Z) \oplus H_1(X,\Z)$;

\item $\pi_2(\mathcal{N}^s_\alpha(n_1,n_2,d_1,d_2))$ is the middle term in an exact sequence
\begin{equation*}
0\lra\Z\lra\pi_2(\mathcal{N}^s_\alpha(n_1,n_2,d_1,d_2))\lra 
Q(n_1,n_2,d_1,d_2)\lra 0\
\end{equation*}

\noindent where

\begin{equation*}
Q(n_1,n_2,d_1,d_2)=
\begin{cases} \Z \oplus \Z \oplus \Z_{\GCD(n_2,d_2)}
\oplus\Z_{\GCD(n_1-n_2,d_1-d_2)}&\text{if $n_2>1$\ and $(n_1-n_2)>1$}\\
\Z \oplus \Z_{\GCD(n_2,d_2)}&\text{if $n_2>1$\ and $(n_1-n_2)=1$}\\
\Z\oplus\Z_{\GCD(n_1-n_2,d_1-d_2)}&\text{if $n_2=1$\ and $(n_1-n_2)>1$}\\
0&\text{if $n_2=1$\ and $n_1=2$}\\
\end{cases}
\end{equation*}
\end{enumerate}
\end{theorem}

\begin{remark}\label{rem:coprime2}
  Theorem~\ref{thm:homotopy-triples} gives a complete description of
  $\pi_2(\mathcal{N}^s_\alpha(n_1,n_2,d_1,d_2))$ when the co-primality
  conditions $\GCD(n_2,d_2) = 1$ and $\GCD(n_1-n_2,n_2-d_2) = 1$ hold,
  and of $\pi_2(\mathcal{N}^s_\alpha(n_1,n_2,d_1,d_2)) \otimes \Q$
  under all circumstances (cf.\ Remark~\ref{rem:coprime-extension}).

  Using the results of \cite{BD}, a complete description of
  $\pi_2(\mathcal{N}^s_\alpha(n_1,n_2,d_1,d_2))$ can also be given in
  the case when $n_2=1$ and $n_1-n_2 > 1$, as we now explain.  In that
  paper, the moduli space of \emph{stable pairs} $(V,\phi)$ was
  studied.  Here $V$ is a vector bundle and $\phi\in H^0(X,V)$ is a
  holomorphic section of $V$. Viewing the section $\phi$ as a map
  $\phi\colon \mathcal{O} \to V$, a pair $(V,\phi)$ gives rise to a
  triple $(E_1,E_2,\phi) = (V,\mathcal{O},\phi)$.  Through this
  correspondence it can be seen that the moduli space of triples
  $\mathcal{N}^s_\alpha$ of triples with $n_2=1$ fibres over the
  Jacobian variety of the curve, with fibres isomorphic to the moduli
  space of pairs.  Among other things, in \cite{BD} the second
  homotopy group of the moduli space of pairs was calculated to be $\Z
  \oplus \Z$ for $\alpha$ between $\alpha_m$ and the first critical
  value of $\alpha$ larger than $\alpha_m$.  It follows from these
  results that, when $d_2=1$, one has $\pi_2(\mathcal{N}^s_\alpha) =
  \Z \oplus \Z $ for such $\alpha$.  Combining this
  fact with Proposition~\ref{homotopy-large-alpha} it follows that
  $\pi_2(\mathcal{N}^s_\alpha) = \Z \oplus \Z $ for $2g-2 \leq \alpha
  < \alpha_M$ in the case $n_2=1$ and $n_1-n_2 > 1$.
\end{remark}

\subsection{Homotopy groups of $\mathcal{M}(p,q, a,b)$}

Combining Propositions \ref{homotopy-minima}, \ref{index-bound} and
\ref{minima-upq} we have the following.

\begin{theorem}\label{homotopy-minima-upq}
Let $\GCD(p+q,a+b)=1$. Then
$$
\pi_i(\mathcal{M}(p,q,a,b))\cong \pi_i(\mathcal{N}(p,q,a,b)),
\;\;\;{for}\;\;\; i\leq 2g-4.
$$
\end{theorem}

\bigskip

As a corollary of Theorems  \ref{homotopy-minima-upq},
\ref{thm:minima=triple-moduli} and \ref{thm:homotopy-triples}
and Proposition \ref{triples-duality},
we conclude  the following.

\begin{theorem}\label{homotopy-higgs-p-q-a-b}
Let $p\neq q$ and $\GCD(p+q,a+b)=1$ and let $g\geq 3$. Then
\begin{enumerate}
\item $\pi_1(\mathcal{M}(p,q,a,b))\cong  H_1(X,\Z) \oplus H_1(X,\Z)$;

\item $\pi_2(\mathcal{M}(p,q,a,b))$ is the middle term in an exact sequence
\begin{equation*}
0\lra\Z\lra\pi_2(\mathcal{M}(p,q,a,b))\lra Q(n_1,n_2,d_1,d_2)\lra 0\ ,
\end{equation*}

\noindent where

\begin{equation*}
(n_1,n_2,d_1,d_2)=
\begin{cases} (p,q,a+p(2g-2),b)&\text{if $\tau<0$\ and $p>q$}\\
(q,p,-b,-a-p(2g-2))&\text{if $\tau<0$\ and $p<q$}\\
(p,q,-a,-b-q(2g-2))&\text{if $\tau>0$\ and $p>q$}\\
(q,p,b+q(2g-2),a)&\text{if $\tau>0$\ and $p<q$}\\
\end{cases}
\end{equation*}

\noindent and where $Q(n_1,n_2,d_1,d_2)$ is as in Theorem
\ref{thm:homotopy-triples}.
\end{enumerate}
\end{theorem}

\begin{remark}
  Theorem~\ref{homotopy-higgs-p-q-a-b} gives a complete description of
  $\pi_2(\mathcal{M}(p,q,a,b))$ when the co-primality
  conditions $\GCD(n_2,d_2) = 1$ and $\GCD(n_1-n_2,n_2-d_2) = 1$ hold,
  and of $\pi_2(\mathcal{M}(p,q,a,b))$ $\otimes \Q$
  under all circumstances (cf.\ Remarks~\ref{rem:coprime-extension}
  and \ref{rem:coprime2}).
\end{remark}

\begin{remark}
As a consequence of Theorem \ref{homotopy-minima-upq} and the connectedness
of $\mathcal{N}(p,q,a,b)$ we have that
$\mathcal{M}(p,q,a,b)$ is also connected,  as proved in \cite{BGG1}.
\end{remark}


\end{document}